\documentclass{article}
\usepackage{latexsym}
\usepackage{amsmath}
\usepackage{amssymb}
\usepackage{amsthm}
\usepackage{authblk}
\newtheorem{theorem}{Theorem}[section]
\newtheorem{corollary}[theorem]{Corollary}
\newtheorem{lemma}[theorem]{Lemma}

\theoremstyle{remark}

\theoremstyle{definition}
\newtheorem{definition}[theorem]{Definition}

\newcommand{\Symm}[1]{\mathcal{S}^{#1}}
\newcommand{\copos}[1]{\mathcal{COP}^{#1}}

\newcommand{\MSOZ}[1]{\mathcal{V}^{#1}_{\min}}

\DeclareMathOperator\Supp{supp}
\DeclareMathOperator\Esupp{esupp}
\DeclareMathOperator\Comp{comp}

\begin{document}

\title{On the structure of the $6 \times 6$ copositive cone}

\author[1]{\small Roland Hildebrand}
\author[2]{\small Andrey Afonin}

\affil[1]{\footnotesize Univ.~Grenoble Alpes, CNRS, Grenoble INP, LJK, 38000 Grenoble, France}
\affil[2]{\footnotesize EPFL, Systems communications section, 1015 Lausanne, Switzerland}

\maketitle

\begin{abstract}
In this work we complement the description of the extreme rays of the $6 \times 6$ copositive cone with some topological structure. In a previous paper we decomposed the set of extreme elements of this cone into a disjoint union of pieces of algebraic varieties of different dimension. In this paper we link this classification to the recently introduced combinatorial characteristic called extended minimal zero support set. We determine those components which are essential, i.e., which are not embedded in the boundary of other components. This allows to drastically decrease the number of cases one has to consider when investigating different properties of the $6 \times 6$ copositive cone. As an application, we construct an example of a copositive $6 \times 6$ matrix with unit diagonal which does not belong to the Parrilo inner sum of squares relaxation ${\cal K}^{(1)}_6$.
\end{abstract}



\section{Introduction}

An element $A$ of the space $\Symm{n}$ of real symmetric $n \times n$ matrices is called \emph{copositive} if $x^TAx \geq 0$ for all vectors $x \in \mathbb R_+^n$. The set of such matrices forms the \emph{copositive cone} $\copos{n}$. This cone plays an important role in non-convex optimization, as many difficult optimization problems can be reformulated as conic programs over $\copos{n}$. For a detailed survey of the applications of this cone see, e.g.,~\cite{Duer10,Bomze12,BSU12,HUS10,DuerRendl21}.

An important characteristic of the copositive cone in relation to optimization over $\copos{n}$ are its extreme rays. Knowledge of the extreme rays allows, e.g., to check the exactness of tractable inner relaxations of the cone~\cite{DDGH13b}. While the extreme rays of $\copos{n}$ for $n \leq 5$ are well understood~\cite{HallNewman63},\cite{Hildebrand12a}, those of $\copos{6}$ are completely classified~\cite{AfoninHildebrandDickinson21}, but the geometry and importance of different types of extreme elements is still obscure. The classification consists of a finite list of parameterized analytic expressions together with constraints on the parameters. Thus the extreme elements of $\copos{6}$ constitute a finite union of pieces of analytic (in fact, algebraic) manifolds. 

In~\cite{AfoninHildebrandDickinson21} the main tool for the classification of the extreme elements of $\copos{6}$ was the \emph{minimal zero support set} introduced in~\cite{Hildebrand14a}. However, it turned out that there can be several types of extreme elements, described by different parameterized analytic expressions, which nevertheless share the same minimal zero support set. Hence this characteristic is too coarse to explain the decomposition. A more adequate characteristic has been developed in~\cite{Hildebrand20a}, the \emph{extended minimal zero support set}. The latter was implicitly discovered already in~\cite{DickinsonHildebrand16}. Elements of $\copos{n}$ which share the same extended minimal zero support set are described by the same polynomial relations on their matrix entries, and hence lie on one algebraic manifold and are described by the same parameterized analytic formula.

In this paper we review the classification of the exceptional extreme rays of $\copos{6}$ in connection to the extended minimal zero support set. It turns out that one given piece can contain matrices with different extended minimal zero support sets, so the correspondence between pieces and support sets is not one-to-one. However, for every piece there exists a single (we shall call it \emph{main}) extended minimal zero support set which characterizes almost all matrices in this piece, and the remaining matrices with different support sets lie on the boundary of the piece. In Section \ref{sec:classification} we provide this main support set for each of the pieces from~\cite{AfoninHildebrandDickinson21}. The results are summarized in Table~\ref{tab1}.

For some purposes, when using the classification, it is sufficient to consider only the main extended zero support sets. Suppose, e.g., that we want to check inclusion of the extreme elements of $\copos{6}$ into some closed set $C$, and this is most simply done by examinating the extended minimal zero support sets of extreme matrices. Then it is sufficient to check only the main support sets, because inclusion of the relative interior of a piece in $C$ implies also inclusion of its boundary by closedness of $C$. 

However, this idea can be carried further, because one piece of extreme matrices can lie entirely on the boundary of another piece, and in this case we do not need to consider the former. One of the results presented in this paper is to single out those pieces of extreme matrices of $\copos{6}$ which do not lie on the boundary of another such piece. We shall call such a piece \emph{essential}. It turns out that only few of the main components are essential, which for some purposes leads to a significant reduction in complexity. The main tool for detecting the considered topological relation between the components is the extended minimal zero support set. The list of essential components is deduced in Section~\ref{sec:c6} and summarized in Theorem~\ref{th:essential}. It may also be helpful in the study of the $6 \times 6$ completely positive cone~\cite{ShakedMonderer17}.

We apply the obtained reduction in complexity to investigate the exactness of the Parrilo ${\cal K}_6^{(1)}$ inner approximating cone (see~\cite{ParriloThesis}) on the unit diagonal affine section of $\copos{6}$. It is known that for $n \leq 5$ the corresponding affine sections of $\copos{n}$ and ${\cal K}_n^{(1)}$ coincide~\cite{DDGH13b}, and that they do not coincide for $n \geq 7$~\cite{LaurentVargas22}. In Section~\ref{sec:ParriloK1} we give a negative answer for the remaining case $n = 6$ by presenting an example of a matrix in $\copos{6}$ with unit diagonal which does not lie in ${\cal K}_6^{(1)}$ (Theorem \ref{th:counterexample}).

\subsection{Notations and preliminaries} \label{subs:notation}

The space of real symmetric matrices of size $n \times n$ will be denoted by $\Symm{n}$, the cone of positive semi-definite such matrices by ${\cal S}_+^n$.

A copositive matrix $A \in \copos{n}$ which is not the sum of a positive semi-definite matrix and an entry-wise nonnegative matrix is called \emph{exceptional}. A non-zero copositive matrix $A \in \copos{n}$ is called \emph{extremal} if whenever $A = A_1 + A_2$ with $A_1,A_2 \in \copos{n}$, the summands $A_1,A_2$ must be nonnegative multiples of $A$.

For an index set $I \subset \{1,\dots,n\}$, denote by $\overline{I}$ its complement $\{1,\dots,n\} \setminus I$.

We shall denote vectors with lower-case letters and matrices with upper-case letters. Individual entries of a vector $u$ and a matrix $A$ will be denoted by $u_i$ and $A_{ij}$ respectively. For a matrix $A$ and a vector $u$ of compatible dimension, the $i$-th element of the matrix-vector product $Au$ will be denoted by $(Au)_i$. Inequalities $u \geq \bf0$ on vectors will be meant element-wise, where we denote by ${\bf 0} = (0,\dots,0)^T$ the all-zeros vector. Similarly we denote by ${\bf 1} = (1,\dots,1)^T$ the all-ones vector. We further let $e_i$ be the unit vector with $i$-th entry equal to one and all other entries equal to zero. For a subset $I \subset \{1,\dots,n\}$ we denote by $A_{I}$ the principal submatrix of $A$ whose elements have row and column indices in $I$, i.e. $A_{I} = (A_{ij})_{i,j \in I}\in\Symm{|I|}$. For subsets $I,J \subset \{1,\dots,n\}$ we denote by $A_{I \times J}$ the submatrix of $A$ whose elements have row indices in $I$ and column indices in $J$. Similarly for a vector $u \in \mathbb R^n$ we define the subvector $u_{I} = (u_i)_{i \in I}\in\mathbb R^{|I|}$.

Let $\Delta = \{ u \in \mathbb R_+^n \mid {\bf 1}^Tu = 1 \}$ be the standard simplex.

For a nonnegative vector $u \in \mathbb R_+^n$ we define its {\it support} as $\Supp{u} = \{i\in\{1,\ldots,n\}\mid u_i > 0 \}$.

A zero $u$ of a copositive matrix $A$ is called {\it minimal} if there exists no zero $v$ of $A$ such that the inclusion $\Supp{v} \subset \Supp{u}$ holds strictly. We shall denote the set of minimal zeros of a copositive matrix $A$ by $\MSOZ{A}$ and the ensemble of supports of the minimal zeros of $A$ by $\Supp{\MSOZ{A}}$. To each index set $I$ there exists at most one minimal zero $u \in \Delta$ of $A$ with $\Supp{u} = I$ \cite[Lemma 3.5]{Hildebrand14a}, hence the minimal zero support set $\Supp{\MSOZ{A}}$ is in bijective correspondence to the minimal zeros of $A$ which are contained in $\Delta$.

For a zero $u$ of a copositive matrix $A$, the matrix-vector product $Au$ is nonnegative. We call the set $\Comp{u} := \overline{\Supp(Au)}$ the \emph{complementary index set} of the zero, and the pair $\Esupp{u} = (\Supp{u},\Comp{u})$ of index sets the \emph{extended support}. The ensemble of extended supports of the minimal zeros of $A$ will be called the \emph{extended minimal zero support set} and denoted by $\Esupp\MSOZ{A}$.

By \cite[Lemma 2.5]{DDGH13a} we have that $\Supp{u} \subset \Comp{u}$ for every zero $u$ of a copositive matrix $A$.

We now briefly recollect the necessary results from \cite{Hildebrand20a}. Let ${\cal E} = \{(I_{\alpha},J_{\alpha})\}_{{\alpha} = 1,\dots,m}$ be a collection of pairs of index sets. Define the sets
\[ S_{\cal E} = \{ A \in \copos{n} \mid \Esupp{\MSOZ{A}} = {\cal E} \}, \qquad
Z_{\cal E} = \{ A \in \Symm{n} \mid A_{I_{\alpha} \times J_{\alpha}}\mbox{ is rank deficient }\forall\ {\alpha} = 1,\dots,m \}.
\]
The set $Z_{\cal E}$ is algebraic, given by the zero locus of a finite number of determinantal polynomials. The set $S_{\cal E}$ is a relatively open subset of $Z_{\cal E}$ \cite[Corollary 1]{Hildebrand20a}.

If $Z_{\cal E}$ is irreducible, then either $S_{\cal E}$ does not contain extremal matrices at all, or all matrices in $S_{\cal E}$ are extremal except possibly in an algebraic subset \cite{Hildebrand20a}. Hence the set of extremal matrices of $\copos{n}$ is a finite union of such sets $S_{\cal E}$, possibly minus some submanifolds of lower dimension. We shall call the sets $S_{\cal E}$ in this decomposition \emph{components} of the set of extremal matrices of $\copos{n}$. Extremal matrices in the same component $S_{\cal E}$ share the extended minimal zero support set ${\cal E}$ and hence many properties in connection to zeros and positive semi-definiteness of submatrices. 

Some components $S_{\cal E}$ may be contained in the closure of others. This defines a hierarchical structure on the ensemble of these components and motivates the following notion \cite{AfoninHildebrandDickinson21}.

{\definition \label{def:essential} A component of extremal matrices of $\copos{n}$ is called \emph{essential} if it is not contained in the closure of another such component. }

This notion is useful in the following situation.

{\lemma Let $C \subset {\cal S}^n$ be a closed convex set. Then $\copos{n} \subset C$ if and only if all essential components of extremal matrices of $\copos{n}$ are contained in $C$. }

\begin{proof}
Clearly if $\copos{n} \subset C$, then all extremal matrices of $\copos{n}$ are contained in $C$, and hence all components of extremal matrices.

Let us show the reverse implication. Since $C$ is closed, it contains the boundaries of all essential components. But then by definition all non-essential components of extreme matrices are contained in $C$, and $C$ contains all extreme matrices of $\copos{n}$. The desired inclusion then follows from convexity of $C$.
\end{proof}

\section{Extended minimal zero support sets of main components} \label{sec:classification}

In this section we determine the generic extended minimal zero support set for the pieces of exceptional extreme matrices of $\copos{6}$ obtained in~\cite{AfoninHildebrandDickinson21}. These sets are provided in Table \ref{tab1}.

The extreme matrices of $\copos{6}$ have been classified in~\cite{AfoninHildebrandDickinson21} with respect to their minimal zero support sets, which are listed in Table 1 of that paper. For each case, the corresponding extreme matrices are given by one or several analytic expressions which are parameterized by the diagonal elements of some scaling matrix $D$ and may be by some angles $\phi_j$ (\cite[Section 5]{AfoninHildebrandDickinson21}). The minimal zeros of these matrices depend on the same parameters.

However, matrices having the same minimal zero support set can belong to different components $S_{\cal E}$, i.e., have different extended minimal zero support sets. To determine the extended minimal zero support set of an analytically given copositive matrix $A$ with known minimal zeros is straightforward. To this end, one has to determine the complement $J_{\alpha}$ of the support of the product $Au^{\alpha}$ for each minimal zero $u^{\alpha}$ of $A$, i.e., to find those indices $i$ such that $(Au^{\alpha})_i = 0$. The expression $(Au^{\alpha})_i$ depends on the parameters entering in the analytic expression of $A$ and $u^{\alpha}$. In particular, whether $(Au^{\alpha})_i = 0$ depends only on the angles $\phi_j$. Therefore different values of the parameters $\phi_j$ may yield different index sets $J_{\alpha}$ and hence a different extended minimal zero support set $\Esupp\MSOZ{A}$. Larger complementary supports correspond to additional equality relations on $\phi_j$, and the corresponding component of extreme matrices has lower dimension.

The derivation of the support sets is straightforward, and we shall demonstrate it in detail on one example, namely the extreme matrices corresponding to minimal zero support set 13 in~\cite[Table 1]{AfoninHildebrandDickinson21}.

\bigskip

{\bf Example:} Minimal zero support set $I = \{\{1,2,3\},\{2,3,4\},\{3,4,5\},\{4,5,6\},\{1,5,6\},\{1,2,6\}\}$. The extreme matrices with this support set are of the form $X = DAD$, where
\begin{equation} \label{A13_1} 
A = \begin{pmatrix}
    \scriptstyle 1 & \scriptstyle -\cos\phi_1 & \scriptstyle \cos(\phi_1+\phi_2) & \scriptstyle -\cos(\phi_1+\phi_2+\phi_3) & \scriptstyle \cos(\phi_5+\phi_6) & \scriptstyle -\cos\phi_6 \\
    \scriptstyle -\cos\phi_1 & \scriptstyle 1 & \scriptstyle -\cos\phi_2 & \scriptstyle \cos(\phi_2+\phi_3) & \scriptstyle  A_{25} & \scriptstyle  \cos(\phi_1+\phi_6) \\
    \scriptstyle \cos(\phi_1+\phi_2) & \scriptstyle  -\cos\phi_2 & \scriptstyle  1 & \scriptstyle  -\cos\phi_3 & \scriptstyle  \cos(\phi_3+\phi_4) & \scriptstyle  -\cos(\phi_3+\phi_4+\phi_5) \\
    \scriptstyle -\cos(\phi_1+\phi_2+\phi_3) & \scriptstyle  \cos(\phi_2+\phi_3) & \scriptstyle  -\cos\phi_3 & \scriptstyle  1 & \scriptstyle  -\cos\phi_4 & \scriptstyle  \cos(\phi_4+\phi_5) \\
    \scriptstyle \cos(\phi_5+\phi_6) & \scriptstyle  A_{25} & \scriptstyle  \cos(\phi_3+\phi_4) & \scriptstyle  -\cos\phi_4 & \scriptstyle  1 & \scriptstyle  -\cos\phi_5 \\
    \scriptstyle -\cos\phi_6 & \scriptstyle  \cos(\phi_1+\phi_6) & \scriptstyle  -\cos(\phi_3+\phi_4+\phi_5) & \scriptstyle  \cos(\phi_4+\phi_5) & \scriptstyle  -\cos\phi_5 & \scriptstyle  1
\end{pmatrix},
\end{equation}
$D$ is a positive definite diagonal matrix, $\phi_i > 0$, $\sum_{j=1}^6 \phi_j < 2\pi$, $\phi_i+\phi_{i+1} < \pi$, $i = 1,\dots,5$, $\phi_1+\phi_6 < \pi$, $\phi_1+\phi_2+\phi_3 \geq \phi_4+\phi_5+\phi_6$, $\phi_3+\phi_4+\phi_5 \geq \phi_1+\phi_2+\phi_6$. In addition, if $\phi_2+\phi_3+\phi_4 \ge \phi_1 + \phi_5 + \phi_6$, then $A_{25} = -\cos(\phi_2+\phi_3+\phi_4)$, and if $\phi_2+\phi_3+\phi_4 \le \phi_1 + \phi_5 + \phi_6$, then $A_{25} = -\cos(\phi_1+\phi_5+\phi_6)$. Furthermore, $\sum_{j=1}^6 \phi_j \not= \pi$, or at least two of the three non-strict inequalities are equalities. The two cases with different expression for $A_{25}$ have been listed in~\cite{AfoninHildebrandDickinson21} as 13.1 and 13.2, respectively.

If $u$ is a minimal zero of $A$, then $D^{-1}u$ is a minimal zero of $X$, and the supports of $u$ and $D^{-1}u$ as well as $Au$ and $XD^{-1}u = DAu$ are the same. Therefore the extended minimal zero support set of $X$ equals that of $A$. The minimal zeros $u^1,\dots,u^6$ of $A$ are given by the columns of the matrix
\begin{equation} \label{minimalZeros}
U = \begin{pmatrix} \sin\phi_2 & 0 & 0 & 0 & \sin\phi_5 & \sin(\phi_1+\phi_6) \\ \sin(\phi_1+\phi_2) & \sin\phi_3 & 0 & 0 & 0 & \sin\phi_6 \\ \sin\phi_1 & \sin(\phi_2+\phi_3) & \sin\phi_4 & 0 & 0 & 0 \\ 0 & \sin\phi_2 & \sin(\phi_3+\phi_4) & \sin\phi_5 & 0 & 0 \\ 0 & 0 & \sin\phi_3 & \sin(\phi_4+\phi_5) & \sin\phi_6 & 0 \\ 0 & 0 & 0 & \sin\phi_4 & \sin(\phi_5+\phi_6) & \sin\phi_1 \end{pmatrix}.
\end{equation}
Let us compute the extended support sets of these zeros. The products $Au^{\alpha}$ appear in the columns of the matrix product $AU$, which equals
\[ \begin{pmatrix}  0 &  0 & c_{56}\sin(\phi_3) & c_{456}\sin(\phi_5) &    0 &   0 \\
  0 &  0 &  0 & c_{16}\sin(\phi_4) & c_{156}\sin(\phi_6) &    0 \\
  0 &   0 &   0 &   0 & c_{12}\sin(\phi_5) &  c_{126}\sin(\phi_1) \\
  0 &   0 &   0 &  0 & c_{456}\sin(\phi_5) & c_{45}\sin(\phi_1) \\
 c_{56}\sin(\phi_2) &  0 &   0 &  0 &   0 & c_{156}\sin(\phi_6) \\
 c_{126}\sin(\phi_1) & c_{16}\sin(\phi_3) &    0 &   0 &    0 &   0 \end{pmatrix}
\]
in case 13.1 and
\[ \begin{pmatrix}  0 &  0 & c_{56}\sin(\phi_3) & c_{456}\sin(\phi_5) &    0 &   0 \\
  0 &  0 &  c_{234}\sin(\phi_3) & c_{23}\sin(\phi_5) & 0 &    0 \\
  0 &   0 &   0 &   0 & c_{12}\sin(\phi_5) &  c_{126}\sin(\phi_1) \\
  0 &   0 &   0 &  0 & c_{456}\sin(\phi_5) & c_{45}\sin(\phi_1) \\
 c_{34}\sin(\phi_1) & c_{234}\sin(\phi_3) &  0 &   0 &  0 &   0 \\
 c_{126}\sin(\phi_1) & c_{16}\sin(\phi_3) &    0 &   0 &    0 &   0 \end{pmatrix}
\]
in case 13.2. Here we denoted for brevity
\[ c_{ij} = \cos(\phi_i + \phi_j) - \cos\sum_{k \not\in \{i,j\}} \phi_k, \quad c_{ijk} = \cos(\phi_i + \phi_j + \phi_k) - \cos\sum_{l \not\in \{i,j,k\}}\phi_l.
\]

Note that $\sin\phi_i > 0$ for all $i$, and the expressions $c_{ij}$ are also positive wherever they appear, because if $0 < \phi < \phi'$ with $\phi + \phi' < 2\pi$, then $\cos\phi > \cos\phi'$. Hence the complementary supports can change only in dependence on whether the expressions $c_{ijk}$ are zero or positive. We obtain for the complementary supports
\[ J_1 \setminus I_1 = \left\{ \begin{array}{rcl} \{ 4 \}, & \ & \phi_3+\phi_4+\phi_5 > \phi_1+\phi_2+\phi_6, \\ \{ 4,6 \}, & \ & \phi_3+\phi_4+\phi_5 = \phi_1+\phi_2+\phi_6, \end{array} \right. \quad
J_2 \setminus I_2 = \left\{ \begin{array}{rcl} \{ 1,5 \}, & \ & \phi_2+\phi_3+\phi_4 \ge \phi_1 + \phi_5 + \phi_6, \\ \{ 1 \}, & \ & \phi_2+\phi_3+\phi_4 < \phi_1 + \phi_5 + \phi_6, \end{array} \right.
\]
\[ J_3 \setminus I_3 = \left\{ \begin{array}{rcl} \{ 6 \}, & \ & \phi_2+\phi_3+\phi_4 < \phi_1 + \phi_5 + \phi_6, \\ \{ 2,6 \}, & \ & \phi_2+\phi_3+\phi_4 \geq \phi_1 + \phi_5 + \phi_6, \end{array} \right. \quad
J_4 \setminus I_4 = \left\{ \begin{array}{rcl} \{ 1,3 \}, & \ & \phi_1+\phi_2+\phi_3 = \phi_4+\phi_5+\phi_6, \\ \{ 3 \}, & \ & \phi_1+\phi_2+\phi_3 > \phi_4+\phi_5+\phi_6, \end{array} \right.
\]
\[ J_5 \setminus I_5 = \left\{ \begin{array}{rcl} \{ 4 \}, & \ & \phi_2+\phi_3+\phi_4 > \phi_1 + \phi_5 + \phi_6,\ \phi_1+\phi_2+\phi_3 = \phi_4+\phi_5+\phi_6, \\ \{ 2,4 \}, & \ & \phi_2+\phi_3+\phi_4 \le \phi_1 + \phi_5 + \phi_6,\ \phi_1+\phi_2+\phi_3 = \phi_4+\phi_5+\phi_6, \\ \emptyset, & \ & \phi_2+\phi_3+\phi_4 > \phi_1 + \phi_5 + \phi_6,\ \phi_1+\phi_2+\phi_3 > \phi_4+\phi_5+\phi_6, \\ \{ 2 \}, & \ & \phi_2+\phi_3+\phi_4 \le \phi_1 + \phi_5 + \phi_6,\ \phi_1+\phi_2+\phi_3 > \phi_4+\phi_5+\phi_6,\end{array} \right. 
\]
\[ J_6 \setminus I_6 = \left\{ \begin{array}{rcl} \{ 3 \}, & \ & \phi_2+\phi_3+\phi_4 > \phi_1 + \phi_5 + \phi_6,\ \phi_3+\phi_4+\phi_5 = \phi_1+\phi_2+\phi_6, \\ \{ 3,5 \}, & \ & \phi_2+\phi_3+\phi_4 \le \phi_1 + \phi_5 + \phi_6,\ \phi_3+\phi_4+\phi_5 = \phi_1+\phi_2+\phi_6, \\ \emptyset, & \ & \phi_2+\phi_3+\phi_4 > \phi_1 + \phi_5 + \phi_6,\ \phi_3+\phi_4+\phi_5 > \phi_1+\phi_2+\phi_6, \\ \{ 5 \}, & \ & \phi_2+\phi_3+\phi_4 \le \phi_1 + \phi_5 + \phi_6,\ \phi_3+\phi_4+\phi_5 > \phi_1+\phi_2+\phi_6.\end{array} \right.
\]
Since always $I_{\alpha} \subset J_{\alpha}$, we provided only the differences of the two sets for brevity.

Thus the extended minimal zero support set depends on whether $\phi_2+\phi_3+\phi_4$ is smaller, equal, or greater than $\phi_1 + \phi_5 + \phi_6$, on whether $\phi_3+\phi_4+\phi_5$ is greater or equal to $\phi_1+\phi_2+\phi_6$, and on whether $\phi_1+\phi_2+\phi_3$ is greater or equal to $\phi_4+\phi_5+\phi_6$. Together this yields 12 different combinations, so the set of extreme matrices with minimal zero support set $I$ decomposes into 12 components $S_{\cal E}$ for different extended minimal zero support sets ${\cal E}$. Some of these components are equivalent under permutations of the index set $\{1,\dots,6\}$.

Clearly only the combinations with exclusively strict inequalities on the $\phi_i$ yield components which do not lie on the boundary of components generated by other combinations. We shall call such components \emph{main components}. We hence obtain as the main component for case 13.1 the complementary supports
\[ J_1 = \{1,2,3,4\},\ J_2 = \{1,2,3,4,5\},\ J_3 = \{2,3,4,5,6\},\ J_4 = \{3,4,5,6\},\ J_5 = \{1,5,6\},\ J_6 = \{1,2,6\},
\]
corresponding to the combination $\phi_2+\phi_3+\phi_4 > \phi_1 + \phi_5 + \phi_6$, $\phi_3+\phi_4+\phi_5 > \phi_1+\phi_2+\phi_6$, $\phi_1+\phi_2+\phi_3 > \phi_4+\phi_5+\phi_6$, and for case 13.2
\[ J_1 = \{ 1,2,3,4 \},\ J_2 = \{1,2,3,4 \},\ J_3 = \{ 3,4,5,6 \},\ J_4 = \{ 3,4,5,6 \},\ J_5 = \{ 1,2,5,6 \},\ J_6 = \{ 1,2,5,6 \},
\]
corresponding to the combination $\phi_2+\phi_3+\phi_4 < \phi_1 + \phi_5 + \phi_6$, $\phi_3+\phi_4+\phi_5 > \phi_1+\phi_2+\phi_6$, $\phi_1+\phi_2+\phi_3 > \phi_4+\phi_5+\phi_6$. Thus each of the two pieces are contained in the closure of a single component $S_{\cal E}$ of extreme matrices, which justifies the notation "main".

\bigskip

It turns out that this holds for all pieces of exceptional extreme matrices in the classification in~\cite{AfoninHildebrandDickinson21}: almost all matrices in each piece belong to a single component, while the rest is located on the boundary of this component. As in the above example, one computes the extended minimal zero support sets for these main components of the other pieces. The result is provided in
Table \ref{tab1}. For brevity we present only the differences $J_{\alpha} \setminus I_{\alpha}$, along with the dimension of each component, taken from \cite[Table 2]{AfoninHildebrandDickinson21}.

\begin{table}[t] \centering
{\small
\begin{tabular}{|c|c|l|c|}
\hline
No. & {\small sets} & {\small main extended minimal zero support set} & dim \\
\hline
O5 & $I_{\alpha}$ & \{1,2,3\},\{2,3,4\},\{3,4,5\},\{1,4,5\},\{1,2,5\},\{6\} & 10 \\
& $J_{\alpha} \setminus I_{\alpha}$ & $\emptyset$,$\emptyset$,$\emptyset$,$\emptyset$,$\emptyset$,\{1,2,3,4,5\} &\\
\hline
1 & $I_{\alpha}$ & \{1,2\},\{1,3\},\{1,4\},\{2,5\},\{3,6\},\{4,5,6\} & 8 \\
& $J_{\alpha} \setminus I_{\alpha}$ & \{3,4,5\},\{2,4,6\},\{2,3\},\{1,6\},\{1,5\},$\emptyset$ & \\
\hline
2 & $I_{\alpha}$ & \{1,2\},\{1,3\},\{1,4\},\{2,5\},\{3,5,6\},\{4,5,6\} & 9 \\
& $J_{\alpha} \setminus I_{\alpha}$ & \{3,4,5\},\{2,4,6\},\{2,3\},\{1,6\},$\emptyset$,$\emptyset$ & \\
\hline
3 & $I_{\alpha}$ & \{1,2\},\{1,3\},\{1,4\},\{2,5,6\},\{3,5,6\},\{4,5,6\} & 10 \\
& $J_{\alpha} \setminus I_{\alpha}$ & \{3,4,5\},\{2,4\},\{2,3,6\},$\emptyset$,$\emptyset$,$\emptyset$ & \\
\hline
4 & $I_{\alpha}$ & \{1,2\},\{1,3\},\{2,4\},\{3,4,5\},\{1,5,6\},\{4,5,6\} & 10 \\
& $J_{\alpha} \setminus I_{\alpha}$ & \{3,4,6\},\{2,6\},\{1,5\},$\emptyset$,$\emptyset$,$\emptyset$ & \\
\hline
5 & $I_{\alpha}$ & \{1,2\},\{1,3\},\{1,4,5\},\{2,4,6\},\{3,4,6\},\{4,5,6\} & 11 \\
& $J_{\alpha} \setminus I_{\alpha}$ & \{3,5\},\{2,5,6\},$\emptyset$,$\emptyset$,$\emptyset$,$\emptyset$ & \\
\hline
6 & $I_{\alpha}$ & \{1,2\},\{1,3\},\{2,4,5\},\{3,4,5\},\{2,4,6\},\{3,5,6\} & 11 \\
& $J_{\alpha} \setminus I_{\alpha}$ & \{3,4\},\{2,5,6\},$\emptyset$,$\emptyset$,$\emptyset$,\{1\} & \\
\hline
7 & $I_{\alpha}$ & \{1,5\},\{2,6\},\{1,2,3\},\{2,3,4\},\{3,4,5\},\{4,5,6\} & \\
& $J_{\alpha} \setminus I_{\alpha}$ & \{2,4\},\{1,3\},\{6\},$\emptyset$,$\emptyset$,$\emptyset$ & 11 \\
\hline
8 & $I_{\alpha}$ & \{1,2\},\{1,3,4\},\{1,3,5\},\{2,4,6\},\{3,4,6\},\{2,5,6\} & \\
& $J_{\alpha} \setminus I_{\alpha}$ & \{3,6\},\{5\},\{4\},$\emptyset$,$\emptyset$,$\emptyset$ & 12 \\
\hline
9 & $I_{\alpha}$ & \{1,2\},\{1,3,4\},\{1,3,5\},\{2,4,6\},\{3,4,6\},\{4,5,6\} & \\
9.1 & $J_{\alpha} \setminus I_{\alpha}$ & \{3,6\},$\emptyset$,$\emptyset$,\{5\},$\emptyset$,\{2\} & 12 \\
9.2 & $J_{\alpha} \setminus I_{\alpha}$ & \{3,5,6\},$\emptyset$,\{2\},$\emptyset$,$\emptyset$,$\emptyset$ & 12 \\
\hline
10 & $I_{\alpha}$ & \{1,2\},\{1,3,4\},\{1,3,5\},\{2,4,6\},\{3,5,6\},\{4,5,6\} & 12 \\
& $J_{\alpha} \setminus I_{\alpha}$ & \{3,6\},$\emptyset$,$\emptyset$,\{5\},$\emptyset$,\{2\} & \\
\hline
11 & $I_{\alpha}$ & \{1,2,3\},\{1,2,4\},\{1,2,5\},\{1,3,6\},\{2,4,6\},\{3,4,6\} & 12 \\
& $J_{\alpha} \setminus I_{\alpha}$ & \{5\},\{5\},\{3,4\},$\emptyset$,$\emptyset$,\{5\} & \\
\hline
12 & $I_{\alpha}$ & \{1,2,3\},\{1,2,4\},\{1,2,5\},\{1,3,6\},\{2,4,6\},\{3,5,6\} & 13 \\
& $J_{\alpha} \setminus I_{\alpha}$ & $\emptyset$,\{5\},\{4\},$\emptyset$,$\emptyset$,\{4\} & \\
\hline
13 & $I_{\alpha}$ & \{1,2,3\},\{2,3,4\},\{3,4,5\},\{4,5,6\},\{1,5,6\},\{1,2,6\} & \\
13.1 & $J_{\alpha} \setminus I_{\alpha}$ & \{4\},\{1,5\},\{2,6\},\{3\},$\emptyset$,$\emptyset$ & 12 \\
13.2 & $J_{\alpha} \setminus I_{\alpha}$ & \{4\},\{1\},\{6\},\{3\},\{2\},\{5\} & 12 \\
\hline
14 & $I_{\alpha}$ & \{1,2\},\{1,3\},\{1,4\},\{2,5\},\{4,5\},\{3,6\},\{5,6\} & 6 \\
& $J_{\alpha} \setminus I_{\alpha}$ & \{3,4,5\},\{2,4,6\},\{2,3,5\},\{1,4,6\},\{1,2,6\},\{1,5\},\{2,3,4\} & \\
\hline
15 & $I_{\alpha}$ & \{1,2\},\{1,3,4\},\{1,3,5\},\{1,4,6\},\{2,5,6\},\{3,5,6\},\{4,5,6\} & \\
& $J_{\alpha} \setminus I_{\alpha}$ & \{3,4\},\{2\},$\emptyset$,$\emptyset$,$\emptyset$,$\emptyset$,$\emptyset$ & 12 \\
\hline
16 & $I_{\alpha}$ & \{1,2,3\},\{1,2,4\},\{1,2,5\},\{1,3,6\},\{2,4,6\},\{3,4,6\},\{3,5,6\} & \\
& $J_{\alpha} \setminus I_{\alpha}$ & $\emptyset$,$\emptyset$,$\emptyset$,$\emptyset$,$\emptyset$,\{5\},\{4\} & 13 \\
\hline
17 & $I_{\alpha}$ & \{1,2,3\},\{1,2,4\},\{1,2,5\},\{1,3,6\},\{2,4,6\},\{3,5,6\},\{4,5,6\} & 13 \\
& $J_{\alpha} \setminus I_{\alpha}$ & $\emptyset$,$\emptyset$,$\emptyset$,$\emptyset$,$\emptyset$,\{4\},\{3\} & \\
\hline
18 & $I_{\alpha}$ & \{1,2,3\},\{2,3,4\},\{3,4,5\},\{1,4,5\},\{1,2,5\},\{3,4,6\},\{1,4,6\},\{1,2,6\} & 12 \\
& $J_{\alpha} \setminus I_{\alpha}$ & $\emptyset$,$\emptyset$,\{6\},\{6\},\{6\},\{5\},\{5\},\{5\} & \\
\hline
19 & $I_{\alpha}$ & \{3,4,5\},\{1,4,5\},\{1,2,5\},\{1,2,3\},\{1,5,6\},\{2,3,4,6\} & \\
& $J_{\alpha} \setminus I_{\alpha}$ & $\emptyset$,\{6\},$\emptyset$,$\emptyset$,\{4\},$\emptyset$ & 14 \\
\hline
\end{tabular}
}
\caption{Extended minimal support sets ${\cal E} = \{(I_{\alpha},J_{\alpha})\}_{\alpha = 1,\dots,m}$ and dimensions of main components of exceptional extreme matrices in $\copos{6}$. Since $I_{\alpha} \subset J_{\alpha}$, for brevity only $I_{\alpha}$ and $J_{\alpha} \setminus I_{\alpha}$ are given for each minimal zero. }
\label{tab1}
\end{table}

\section{Essential components of extremal matrices in $\copos{6}$} \label{sec:c6}

In this section we investigate which of the obtained 22 main components of exceptional extreme matrices listed in Table~\ref{tab1} are essential, i.e., do not lie on the boundary of other main components. The main tool is the following result, which is a consequence of \cite[Lemma 6]{Hildebrand20a}.

{\lemma \label{lem:essential_crit} Let a component $S_{{\cal E}'}$ of extremal matrices be contained in the closure of another component $S_{\cal E}$. Let ${\cal E} = \{ (I_{\alpha},J_{\alpha}) \}_{{\alpha} = 1,\dots,m}$, ${\cal E}' = \{ (I'_{\alpha},J'_{\alpha}) \}_{{\alpha} = 1,\dots,m'}$. Then for every ${\alpha} = 1,\dots,m$ there exists ${\alpha}' \in \{1,\dots,m'\}$ such that $I'_{{\alpha}'} \subset I_{\alpha}$, $J_{\alpha} \subset J'_{{\alpha}'}$. \qed }

\medskip

Note that as in the classification in~\cite{AfoninHildebrandDickinson21}, Table~\ref{tab1} lists the extended minimal zero support sets of the main components only up to a permutation of the index set $\{1,\dots,6\}$. This means that every entry in Table~\ref{tab1} stands for potentially up to $6! = 720$ different main components, and each of these components can contain another main component from Table~\ref{tab1} in its boundary. We hence have to allow for this freedom when checking the criterion in Lemma \ref{lem:essential_crit} on pairs of extended minimal zero support sets from Table~\ref{tab1}. Direct verification of the criterion in conjunction with a strict inequality between the dimension of the components (the boundary of $S_{\cal E}$ can contain $S_{{\cal E}'}$ only if $\dim\,S_{\cal E} > \dim\,S_{{\cal E}'}$) yields the set of potential pairs satisfying an inclusion relation which is listed in Table~\ref{tab3}.

\begin{table}[t] \centering
{\small
\begin{tabular}{|c|l|}
\hline
No. & {\small may possibly be in the closure of}  \\
\hline
O5 & 8,16,17 \\
1 & 2,3,4,5,6,7,8,9.1,9.2,10,15,16,17 \\
2 & 3,4,5,6,7,8,9.1,9.2,10,15,16,17 \\
3 & 5,9.2,15,16,17 \\
4 & 5,6,7,8,9.1,9.2,10,15,16,17 \\
5 & 9.2,15,16,17 \\
6 & 8,9.1,9.2,10,16,17 \\
7 & 8,9.1,10,15,16,17 \\
8 & 16 \\
9.1 & 16 \\
9.2 & 16,17 \\
10 & 17 \\
11 & 19 \\
12 & 19 \\
13.1 & \\
13.2 & \\
14 & 1,2,3,4,5,6,7,8,9.1,9.2,10,11,12,13.1,13.2,15,16,17,18 \\
15 & 16,17 \\
16 & \\
17 & \\
18 & 12,16,19 \\
19 & \\
\hline
\end{tabular}
}
\caption{Pairs of main components of exceptional extremal matrices in $\copos{6}$ satisfying both the criterion in Lemma \ref{lem:essential_crit} and a strict inequality on the dimensions. }
\label{tab3}
\end{table}

From Table \ref{tab3} it follows that the main components corresponding to cases 13.1,13.2,16,17,19 are essential. For the remaining 17 components we now show that they lie on the boundary of other components and are hence non-essential, by providing an appropriate permutation of the indices $\{1,\dots,6\}$ and an explicit limit.

%In \ref{tab1} we provide you with cases of exceptional extreme matrices $\copos{6}$ and corresponding minimal support sets $\SI$. In the \ref{tab2} we present the dimensions of the mutually non-isomorphic strata of exceptional extremal matrices with unit diagonal corresponding to the minimal zero support sets 1--19 in \ref{tab1}. The respective maximal dimension equals the number of free parameters in the expressions of corresponding matrix \cite{AfoninHildebrandDickinson21}. The strata of smaller dimension are obtained by letting some of the non-strict inequalities on the parameters be equalities. Removing the restriction that the diagonal elements of the matrix equal 1 increases all dimensions by 6.

%\medskip

%\begin{table}[]
%    \centering
%    \begin{tabular}{|c|c||c|c||c|c||c|c|}
%    \hline
%    Case No. & Dim. & Case No. & Dim. & Case No. & Dim. & Case No. & Dim. \\
%    \hline
%    1 & 2 & 6 & 5 & 11 & 6 & 16 & 7,6,6,5 \\
%    2 & 3 & 7 & 5,4 & 12 & 7 & 17 & 7 \\
%    3 & 4 & 8 & 6,5,5,4 & 13 & 6,6,5,5,4,3 & 18 & 6 \\
%    4 & 4 & 9 & 6,6 & 14 & 0 & 19 & 8,7 \\
%    5 & 5 & 10 & 6 & 15 & 6,5,5 & & \\
%    \hline
%    \end{tabular}
%    \caption{Dimensions of strata of exceptional extremal matrices with unit diagonal}
%    \label{tab2}
%\end{table}

\begin{itemize}
    \item \textbf{Case O5} is in the closure of \textbf{Case 18} when the last row and column tend to zero.

    \item \textbf{Case 14} is in the closure of \textbf{Case 1} after the substitution $\phi_1 \rightarrow 0$, $\phi_2\rightarrow 0$.

    \item \textbf{Case 1} is in the closure of \textbf{Case 2} after the substitution $\phi_1 \rightarrow \phi_2$, $\phi_2 \rightarrow 0$, $\phi_3 \rightarrow \pi-\phi_1-\phi_2$.

    \item \textbf{Case 2} is in the closure of \textbf{Case 3} after the permutation $(124356)$ and the substitution $\phi_1 \rightarrow \phi_1$, $\phi_2 \rightarrow \pi-\phi_1$,$\phi_3 \rightarrow \phi_3$, $\phi_4 \rightarrow \phi_2$.

    \item \textbf{Case 3} is in the closure of \textbf{Case 5} after the permutation $(152346)$ and the substitution $\phi_1 \rightarrow \phi_1$,$\phi_2 \rightarrow (\pi-\phi_1-\phi_2)$,$\phi_3 \rightarrow \phi_4$, $\phi_4 \rightarrow \phi_3$,$\phi_5 \rightarrow 0$.

    \item \textbf{Case 10} is in the closure of \textbf{Case 17} after the permutation $(241536)$ and the substitution $\phi_1 \rightarrow \phi_3$, $\phi_2 \rightarrow \phi_2$, $\phi_3 \rightarrow 0$, $\phi_4 \rightarrow \pi-\phi_1-\phi_2$, $\phi_5 \rightarrow \phi_6$, $\phi_6 \rightarrow \phi_5$, $\phi_7 \rightarrow \phi_4$.

    \item \textbf{Case 11} is in the closure of \textbf{Case 19} after the permutation $(514263)$ and the substitution $\phi_1 \rightarrow \phi_4$, $\phi_2 \rightarrow \phi_1$, $\phi_3 \rightarrow \phi_2$, $\phi_4 \rightarrow \phi_3$, $\phi_5 \rightarrow \phi_5$, $\phi_6 \rightarrow \pi-\phi_2-\phi_6$, $\phi_7 \rightarrow \pi-\phi_6+\phi_3$, $a_{24} \rightarrow \cos(\phi_4+\phi_5)$, $a_{36} \rightarrow b_3$.

    \item \textbf{Case 12} is in the closure of \textbf{Case 19} after the permutation $(152463)$ and the substitution $\phi_1 \rightarrow \phi_6$, $\phi_2 \rightarrow \phi_3$, $\phi_3 \rightarrow \phi_2$, $\phi_4 \rightarrow \phi_1$, $\phi_5 \rightarrow \phi_5$, $\phi_6 \rightarrow \phi_4$, $\phi_7 \rightarrow \pi-\phi_5-\phi_7$, $a_{24} \rightarrow b_1$, $a_{36} \rightarrow -\cos\phi_7$.

    \item Part $0 \leq \phi_6 < \phi_2$ of case 18 is in the closure of case 16 after the permutation $(631254)$ and the substitution $\phi_1 \rightarrow \phi_5$, $\phi_2 \rightarrow \phi_1$, $\phi_3 \rightarrow \phi_2-\phi_6$, $\phi_4 \rightarrow \phi_2$, $\phi_5 \rightarrow \phi_4$, $\phi_6 \rightarrow \phi_3+\phi_6$, $\phi_7 \rightarrow \phi_3$. Using the permutation $(213654)$ instead we obtain the part $-\phi_3 < \phi_6 \leq 0$ of case 18, because this part is obtained from the former part by the permutation $(432156)$. As a result, \textbf{Case 18} is in the closure of \textbf{Case 16}.

    \item \textbf{Case 15} is in the closure of \textbf{Case 16} after the permutation $(654321)$ and the substitution $\phi_1 \rightarrow \phi_4$, $\phi_2 \rightarrow \phi_5$, $\phi_3 \rightarrow \phi_3$, $\phi_4 \rightarrow \pi-\phi_5-\phi_6$, $\phi_5 \rightarrow \phi_1$, $\phi_6 \rightarrow \phi_2$, $\phi_7 \rightarrow 0$.

    \item \textbf{Case 9.1} is in the closure of \textbf{Case 16} after the permutation $(241356)$ and the substitution $\phi_1 \rightarrow \phi_1$, $\phi_2 \rightarrow \phi_2$, $\phi_3 \rightarrow 0$, $\phi_4 \rightarrow \pi-\phi_2-\phi_3$, $\phi_5 \rightarrow \phi_4$, $\phi_6 \rightarrow \phi_5$, $\phi_7 \rightarrow \phi_6$.

    \item \textbf{Case 9.2} is in the closure of \textbf{Case 16} after the permutation $(643152)$ and the substitution $\phi_1 \rightarrow \phi_1$, $\phi_2 \rightarrow \phi_4$, $\phi_3 \rightarrow \phi_5$, $\phi_4 \rightarrow \phi_6$, $\phi_5 \rightarrow \phi_2$, $\phi_6 \rightarrow 0$, $\phi_7 \rightarrow \pi-\phi_2-\phi_3$.

    \item \textbf{Case 8} is in the closure of \textbf{Case 16} after the permutation $(316452)$ and the substitution $\phi_1 \rightarrow 0$, $\phi_2 \rightarrow \phi_5$, $\phi_3 \rightarrow \phi_4$, $\phi_4 \rightarrow \phi_6$, $\phi_5 \rightarrow \phi_2$, $\phi_6 \rightarrow \phi_1$, $\phi_7 \rightarrow \phi_3$.

    \item \textbf{Case 5} is in the closure of \textbf{Case 16} after the permutation $(654132)$ and the substitution $\phi_1 \rightarrow \phi_2$, $\phi_2 \rightarrow \phi_1$, $\phi_3 \rightarrow \phi_3$, $\phi_4 \rightarrow \phi_4$, $\phi_5 \rightarrow \phi_5$, $\phi_6 \rightarrow 0$, $\phi_7 \rightarrow 0$.

    \item \textbf{Case 7} is in the closure of \textbf{Case 16} after the permutation $(463125)$ and the substitution $\phi_1 \rightarrow \phi_3$, $\phi_2 \rightarrow \phi_4$, $\phi_3 \rightarrow 0$, $\phi_4 \rightarrow \phi_5$, $\phi_5 \rightarrow \phi_2$, $\phi_6 \rightarrow \phi_1$, $\phi_7 \rightarrow 0$.

    \item \textbf{Case 6} is in the closure of \textbf{Case 16} after the permutation $(426135)$ and the substitution $\phi_1 \rightarrow \phi_3$, $\phi_2 \rightarrow \phi_2$, $\phi_3 \rightarrow 0$, $\phi_4 \rightarrow \pi-\phi_2-\phi_4$, $\phi_5 \rightarrow \phi_1$, $\phi_6 \rightarrow 0$, $\phi_7 \rightarrow \phi_5$.

    \item \textbf{Case 4} is in the closure of \textbf{Case 16} after the permutation $(645213)$ and the substitution $\phi_1 \rightarrow \phi_3$, $\phi_2 \rightarrow \phi_2$, $\phi_3 \rightarrow 0$, $\phi_4 \rightarrow \phi_1$, $\phi_5 \rightarrow \phi_4$, $\phi_6 \rightarrow 0$, $\phi_7 \rightarrow 0$.
\end{itemize}

We obtain the following result.

{\theorem \label{th:essential} Out of the 22 mutually non-equivalent main components of exceptional extremal matrices in $\copos{6}$ exactly the 5 components 13.1,13.2,16,17,19 in Table \ref{tab1} are essential. \qed }

Let us remark that in the case of the cone $\copos{5}$, there are two mutually non-equivalent components of exceptional extremal matrices, of which one is essential.

\section{Relation of $\copos{6}$ and the Parrilo cone ${\cal K}_6^{(1)}$} \label{sec:ParriloK1}

In \cite[Theorem 5.2]{ParriloThesis} the following inner approximation of the copositive cone $\copos{n}$ was established. The cone ${\cal K}_n^{(1)}$ is defined as the set of matrices $A \in {\cal S}^n$ such that there exist symmetric matrices $\Lambda^1,\dots,\Lambda^n \in {\cal S}^n$ satisfying
\begin{equation} \label{K1condition} 
M^i := A - \Lambda^i \succeq 0, \quad \Lambda^i_{ii} = 0, \quad \Lambda^i_{jj} + 2\Lambda^j_{ij} = 0, \quad m_{ijk} := \Lambda^i_{jk} + \Lambda^j_{ik} + \Lambda^k_{ij} \geq 0
\end{equation}
for all distinct $i,j,k \in \{1,\dots,n\}$. In this section we shall construct an extreme matrix $A \in \copos{6} \setminus {\cal K}_6^{(1)}$ with unit diagonal. By \cite[Theorem 5.2]{ParriloThesis} any such matrix must be exceptional.

If such a counterexample exists, it must lie in one of the components of exceptional extreme matrices. But the cone ${\cal K}_6^{(1)}$ is closed, and so is its affine section consisting of unit diagonal matrices. Hence there must be at least one \emph{essential} component which also contains matrices in the difference $A \in \copos{6} \setminus {\cal K}_6^{(1)}$ with unit diagonal. It is hence sufficient to search for a counterexample in the essential components. The counterexample constructed below was found in the main component corresponding to case 13.1.

We briefly outline our strategy. First we show that if $x$ is a zero of a matrix $A \in {\cal K}_n^{(1)}$, then $x$ must be in the kernel of the matrices $M^i$ with $i \in \Supp\,x$. If $A$ is in the main component corresponding to case 13.1 in Table~\ref{tab1}, then these zeros are sufficiently numerous to determine the $M^i$ completely. We then derive explicit conditions on the angles $\phi_i$ such that $A \in {\cal K}_n^{(1)}$ and construct an example which violates these conditions.

{\lemma \label{lem:kernel} Let $A \in {\cal K}_n^{(1)}$, let $x$ be a zero of $A$, and let $M^i,\Lambda^i$ be as in \eqref{K1condition}. Then $M^ix = 0$ for every $i \in \Supp\,x$. }

\begin{proof}
By definition we have $x^TAx = 0$ and hence
\[ 0 \leq \sum_{i=1}^n x_i \cdot (x^TM^ix) = - \sum_{i=1}^n x_i \cdot (x^T\Lambda^ix) = - \sum_{i,j,k=1}^n \Lambda^i_{jk}x_ix_jx_k \leq 0,
\]
because the totally symmetric part of the tensor $\Lambda^i_{jk}$ is element-wise nonnegative, and $x \geq {\bf 0}$. It follows that $x_i \cdot (x^TM^ix) = 0$ for all $i$, which yields the desired conclusion by virtue of $M^i \succeq 0$.
\end{proof}

{\corollary \label{cor:kernel} Let $x^{\alpha},x^{\beta}$ be minimal zeros of $A \in {\cal K}_n^{(1)}$ with extended supports $(I_{\alpha},J_{\alpha})$, $(I_{\beta},J_{\beta})$, respectively, and suppose that $I_{\alpha} \cup I_{\beta} \subset J_{\alpha} \cap J_{\beta}$. Then with $M^i$ as in the previous lemma we have $M^ix^{\alpha} = 0$ for all $i \in I_{\alpha} \cup I_{\beta}$. }

\begin{proof}
Since $I_{\beta} \subset J_{\alpha}$, we have $(x^{\beta})^TAx^{\alpha} = 0$. But then $(x^{\alpha}+x^{\beta})^TA(x^{\alpha}+x^{\beta}) = 0$, and $x^{\alpha}+x^{\beta}$ is a zero of $A$ with support $I_{\alpha} \cup I_{\beta}$. By the previous lemma we obtain that $x^{\alpha}+x^{\beta}$ is in the kernel of $M^i$ for all $i \in I_{\alpha} \cup I_{\beta}$. However, also by this lemma $x^{\alpha}$ is in the kernel of $M^i$ for all $i \in I_{\alpha}$, and $x^{\beta}$ is in the kernel of $M^i$ for all $i \in I_{\beta}$. The claim of the corollary now follows.
\end{proof}

Let now $A$ be of the form \eqref{A13_1} with $A_{25} = -\cos(\phi_2+\phi_3+\phi_4)$ and $\phi_i$ satisfying the conditions
\begin{equation}  \label{phi_cond}
\begin{aligned}
\phi_i > 0, \quad \sum_{j=1}^6 \phi_j < 2\pi, &\quad \phi_i+\phi_{i+1} < \pi, \quad \phi_1+\phi_6 < \pi, \quad \sum_{j=1}^6 \phi_j \not= \pi, \\
\phi_1+\phi_2+\phi_3 > \phi_4+\phi_5+\phi_6, &\quad \phi_3+\phi_4+\phi_5 > \phi_1+\phi_2+\phi_6, \quad \phi_2+\phi_3+\phi_4 > \phi_1 + \phi_5 + \phi_6,
\end{aligned}
\end{equation}
i.e., $A$ is in the main component corresponding to case 13.1 in Table~\ref{tab1}. 

We shall deduce when the inclusion $A \in {\cal K}_6^{(1)}$ is valid. Let $M^i,m_{ijk}$ be as in \eqref{K1condition}. Recall that the minimal zeros $x^1,\dots,x^6$ of $A$ are given by the columns of \eqref{minimalZeros}, and let $(I_{\alpha},J_{\alpha})$ be the extended support set of $x^{\alpha}$. By virtue of Table~\ref{tab1} we have $I_{\alpha} \cup I_{\beta} \subset J_{\alpha} \cap J_{\beta}$ for $(\alpha,\beta) = (1,2),(2,3),(3,4)$. By Lemma \ref{lem:kernel} and Corollary \ref{cor:kernel} we then have that $M^1$ is orthogonal to $x^1,x^2,x^5,x^6$, $M^2$ to $x^1,x^2,x^3,x^6$, $M^3$ and $M^4$ to $x^1,x^2,x^3,x^4$, $M^5$ to $x^2,x^3,x^4,x^5$, and $M^6$ to $x^3,x^4,x^5,x^6$.

It is easily seen that the listed quadruples of minimal zeros are linearly independent, and the $M^i$ have to be of rank at most 2. Direct verification shows that the ranges of $M^i$ are contained the column spaces of the matrices
\[ R^1 = \begin{pmatrix} 1 & 0 \\ -\cos\phi_1 & \sin\phi_1 \\ \cos(\phi_1+\phi_2) & -\sin(\phi_1+\phi_2) \\ -\cos(\phi_1+\phi_2+\phi_3) & \sin(\phi_1+\phi_2+\phi_3) \\ \cos(\phi_5+\phi_6) & \sin(\phi_5+\phi_6) \\ -\cos(\phi_6) & -\sin(\phi_6) \end{pmatrix},\quad R^2 = \begin{pmatrix} -\cos(\phi_1) & -\sin(\phi_1) \\ 1 & 0 \\ -\cos\phi_2 & \sin\phi_2 \\ \cos(\phi_2+\phi_3) & -\sin(\phi_2+\phi_3) \\ -\cos(\phi_2+\phi_3+\phi_4) & \sin(\phi_2+\phi_3+\phi_4) \\ \cos(\phi_1+\phi_6) & \sin(\phi_1+\phi_6) \end{pmatrix},
\]
\[ R^3 = \begin{pmatrix} \cos(\phi_1+\phi_2) & \sin(\phi_1+\phi_2) \\ -\cos(\phi_2) & -\sin(\phi_2) \\ 1 & 0 \\ -\cos\phi_3 & \sin\phi_3 \\ \cos(\phi_3+\phi_4) & -\sin(\phi_3+\phi_4) \\ -\cos(\phi_3+\phi_4+\phi_5) & \sin(\phi_3+\phi_4+\phi_5) \end{pmatrix}, \quad R^4 = \begin{pmatrix} -\cos(\phi_1+\phi_2+\phi_3) & -\sin(\phi_1+\phi_2+\phi_3) \\ \cos(\phi_2+\phi_3) & \sin(\phi_2+\phi_3) \\ -\cos(\phi_3) & -\sin(\phi_3) \\ 1 & 0 \\ -\cos\phi_4 & \sin\phi_4 \\ \cos(\phi_4+\phi_5) & -\sin(\phi_4+\phi_5) \end{pmatrix},
\]
\[ R^5 = \begin{pmatrix} \cos(\phi_5+\phi_6) & -\sin(\phi_5+\phi_6) \\ -\cos(\phi_2+\phi_3+\phi_4) & -\sin(\phi_2+\phi_3+\phi_4) \\ \cos(\phi_3+\phi_4) & \sin(\phi_3+\phi_4) \\ -\cos(\phi_4) & -\sin(\phi_4) \\ 1 & 0 \\ -\cos\phi_5 & \sin\phi_5 \end{pmatrix}, \quad R^6 = \begin{pmatrix} -\cos\phi_6 & \sin\phi_6 \\ \cos(\phi_1+\phi_6) & -\sin(\phi_1+\phi_6) \\ -\cos(\phi_3+\phi_4+\phi_5) & -\sin(\phi_3+\phi_4+\phi_5) \\ \cos(\phi_4+\phi_5) & \sin(\phi_4+\phi_5) \\ -\cos(\phi_5) & -\sin(\phi_5) \\ 1 & 0 \end{pmatrix},
\]
respectively, such that $M^i = R^iD^i(R^i)^T$ with $D^i \in {\cal S}_+^2$. 

The conditions $A_{ii} = 1 = M^i_{ii}$ lead to $D^i_{11} = 1$ for all $i$. Further with
\[ j = \left\{ \begin{array}{rcl} i+1,&\ &i = 1,\dots,5; \\ 1,&& i = 6 \end{array} \right.
\]
the conditions $2A_{ij} = 2M^i_{ij} + 2\Lambda^i_{ij} = 2M^i_{ij} - \Lambda^j_{ii} = 2M^i_{ij} + M^j_{ii} - A_{ii}$ become
\[ -2\cos\phi_i = 2 \cdot \begin{pmatrix} 1 \\ 0 \end{pmatrix}^T \begin{pmatrix} 1 & D^i_{12} \\ D^i_{12} & D^i_{22} \end{pmatrix} \begin{pmatrix} -\cos\phi_i \\ \sin\phi_i \end{pmatrix} + \begin{pmatrix} -\cos\phi_i \\ -\sin\phi_i \end{pmatrix}^T \begin{pmatrix} 1 & D^j_{12} \\ D^j_{12} & D^j_{22} \end{pmatrix} \begin{pmatrix} -\cos\phi_i \\ -\sin\phi_i \end{pmatrix} - 1,
\]
which simplifies to
\[ 2D^i_{12} + 2D^j_{12} \cdot \cos\phi_i + (D^j_{22} - 1) \cdot \sin\phi_i = 0.
\]
Similarly, the conditions $2A_{ij} = 2M^j_{ij} + 2\Lambda^j_{ij} = 2M^j_{ij} - \Lambda^i_{jj} = 2M^j_{ji} + M^i_{jj} - A_{jj}$ become
\[ -2\cos\phi_i = 2 \cdot \begin{pmatrix} 1 \\ 0 \end{pmatrix}^T \begin{pmatrix} 1 & D^j_{12} \\ D^j_{12} & D^j_{22} \end{pmatrix} \begin{pmatrix} -\cos\phi_i \\ -\sin\phi_i \end{pmatrix} + \begin{pmatrix} -\cos\phi_i \\ \sin\phi_i \end{pmatrix}^T \begin{pmatrix} 1 & D^i_{12} \\ D^i_{12} & D^i_{22} \end{pmatrix} \begin{pmatrix} -\cos\phi_i \\ \sin\phi_i \end{pmatrix} - 1,
\]
which simplifies to
\[ -2D^j_{12} - 2D^i_{12} \cdot \cos\phi_i + (D^i_{22} - 1) \cdot \sin\phi_i = 0.
\]
This yields a homogeneous linear system of equations on the quantities $2D^1_{12},\dots,2D^6_{12},D^1_{22}-1,\dots,D^6_{22}-1$ with coefficient matrix
\[ \begin{pmatrix} 1 & \cos\phi_1 & & & & & & \sin\phi_1 & & & & \\ & 1 & \cos\phi_2 & & & & & & \sin\phi_2 & & & \\ & & 1 & \cos\phi_3 & & & & & & \sin\phi_3 & & \\ & & & 1 & \cos\phi_4 & & & & & & \sin\phi_4 & \\ & & & & 1 & \cos\phi_5 & & & & & & \sin\phi_5 \\ \cos\phi_6 & & & & & 1 & \sin\phi_6 & & & & & \\ -\cos\phi_1 & -1 & & & & & \sin\phi_1 & & & & & \\ & -\cos\phi_2 & -1 & & & & & \sin\phi_2 & & & & \\ & & -\cos\phi_3 & -1 & & & & & \sin\phi_3 & & & \\ & & & -\cos\phi_4 & -1 & & & & & \sin\phi_4 & & \\ & & & & -\cos\phi_5 & -1 & & & & & \sin\phi_5 & \\ -1 & & & & & -\cos\phi_6 & & & & & & \sin\phi_6 \end{pmatrix}.
\]
The determinant of this matrix is a non-vanishing trigonometric polynomial in $\phi_1,\dots,\phi_6$, and hence for almost all values of these angles we must have $D^i_{12} = 0$, $D^i_{22} = 1$ for all $i$. This determines the matrices $M^i = R^i(R^i)^T$ and hence the $\Lambda^i = A - R^i(R^i)^T$ completely.

Direct calculation yields that 12 of the 20 independent $m_{ijk}$ are zero, and the remaining are given by
\begin{align*}
\scriptstyle m_{125} &=\, \scriptstyle \cos(\phi_1 + \phi_5 + \phi_6) - \cos(\phi_2 + \phi_3 + \phi_4) + \cos(\phi_5 + \phi_6) - \cos(\phi_1 + \phi_2 + \phi_3 + \phi_4) - \cos(\phi_1) + \cos(\phi_2 + \phi_3 + \phi_4 + \phi_5 + \phi_6), \\
\scriptstyle m_{135} &=\, \scriptstyle [\cos(\phi_1 + \phi_2) - \cos(\phi_3 + \phi_4 + \phi_5 + \phi_6)] + [\cos(\phi_5 + \phi_6) - \cos(\phi_1 + \phi_2 + \phi_3 + \phi_4)] + [\cos(\phi_3 + \phi_4) - \cos(\phi_1 + \phi_2 + \phi_5 + \phi_6)], \\
\scriptstyle m_{136} &=\, \scriptstyle \cos(\phi_1 + \phi_2 + \phi_6) - \cos(\phi_3 + \phi_4 + \phi_5) + \cos(\phi_1 + \phi_2) - \cos(\phi_3 + \phi_4 + \phi_5 + \phi_6) - \cos(\phi_6) + \cos(\phi_1 + \phi_2 + \phi_3 + \phi_4 + \phi_5), \\
\scriptstyle m_{145} &=\, \scriptstyle \cos(\phi_4 + \phi_5 + \phi_6) - \cos(\phi_1 + \phi_2 + \phi_3) + \cos(\phi_5 + \phi_6) - \cos(\phi_1 + \phi_2 + \phi_3 + \phi_4) - \cos(\phi_4) + \cos(\phi_1 + \phi_2 + \phi_3 + \phi_5 + \phi_6), \\
\scriptstyle m_{146} &=\, \scriptstyle \cos(\phi_4 + \phi_5 + \phi_6) - \cos(\phi_1 + \phi_2 + \phi_3) + \cos(\phi_4 + \phi_5) - \cos(\phi_1 + \phi_2 + \phi_3 + \phi_6) - \cos(\phi_6) + \cos(\phi_1 + \phi_2 + \phi_3 + \phi_4 + \phi_5), \\
\scriptstyle m_{236} &=\, \scriptstyle \cos(\phi_1 + \phi_2 + \phi_6) - \cos(\phi_3 + \phi_4 + \phi_5) + \cos(\phi_1 + \phi_6) - \cos(\phi_2 + \phi_3 + \phi_4 + \phi_5) - \cos(\phi_2) + \cos(\phi_1 + \phi_3 + \phi_4 + \phi_5 + \phi_6), \\
\scriptstyle m_{246} &=\, \scriptstyle [\cos(\phi_2 + \phi_3) - \cos(\phi_1 + \phi_4 + \phi_5 + \phi_6)] + [\cos(\phi_1 + \phi_6) - \cos(\phi_2 + \phi_3 + \phi_4 + \phi_5)] + [\cos(\phi_4 + \phi_5) - \cos(\phi_1 + \phi_2 + \phi_3 + \phi_6)], \\
\scriptstyle m_{256} &=\, \scriptstyle \cos(\phi_1 + \phi_5 + \phi_6) - \cos(\phi_2 + \phi_3 + \phi_4) + \cos(\phi_1 + \phi_6) - \cos(\phi_2 + \phi_3 + \phi_4 + \phi_5) - \cos(\phi_5) + \cos(\phi_1 + \phi_2 + \phi_3 + \phi_4 + \phi_6).
\end{align*}

In order for $A$ to be in ${\cal K}_6^{(1)}$ all $m_{ijk}$ have to be nonnegative (and this is actually an equivalence). This is, however, not guaranteed by conditions \eqref{phi_cond}. Note that for $0 < \phi < \phi'$ such that $\phi + \phi' < 2\pi$ we have $\cos\phi - \cos\phi' \geq 0$. This yields nonnegativity of the expressions in square brackets and hence nonnegativity of $m_{135},m_{246}$. The other $m_{ijk}$ may be negative, however, as the following example shows.

For
\[
\phi_1 = 0.20\pi, \quad  \phi_2 = 0.29\pi, \quad \phi_3 = 0.30\pi, \quad \phi_4 = 0.23\pi, \quad \phi_5 = 0.06\pi, \quad \phi_6 = 0.02\pi
\]
we get $m_{136} < -\frac43$, and the above-mentioned determinant of the linear system on $D^i_{12},D^i_{22}$ is non-zero. We obtain the following result.

{\theorem \label{th:counterexample} There exist matrices with unit diagonal in the difference $\copos{6} \setminus {\cal K}_6^{(1)}$. An example is the extreme copositive matrix
\[ \begin{pmatrix} 1 & -\cos(0.20\pi) & \cos(0.49\pi) & -\cos(0.79\pi) & \cos(0.08\pi) & -\cos(0.02\pi) \\
-\cos(0.20\pi) & 1 & -\cos(0.29\pi) & \cos(0.59\pi) & -\cos(0.82\pi) & \cos(0.22\pi) \\
\cos(0.49\pi) & -\cos(0.29\pi) & 1 & -\cos(0.30\pi) & \cos(0.53\pi) & -\cos(0.59\pi) \\
-\cos(0.79\pi) & \cos(0.59\pi) & -\cos(0.30\pi) & 1 & -\cos(0.23\pi) & \cos(0.29\pi) \\
\cos(0.08\pi) & -\cos(0.82\pi) & \cos(0.53\pi) & -\cos(0.23\pi) & 1 & -\cos(0.06\pi) \\
-\cos(0.02\pi) & \cos(0.22\pi) & -\cos(0.59\pi) & \cos(0.29\pi) & -\cos(0.06\pi) & 1 \end{pmatrix}.
\qedhere
\]
}

\section*{Acknowledgements}

The authors are indebted to Peter J.C. Dickinson for his valuable advice in the calculation of the essential extremal strata of the $\copos{6}$ cone.

\bibliographystyle{plain}
\bibliography{copositive,convexity}

\end{document}